# APPLICATIONS OF WALLIS THEOREM

Mihály Bencze[1], Florentin Smarandache[2]

[1]Department of Mathematics, Áprily Lajos College, Braşov, Romania;
[2]Chair of Department of Math & Sciences, University of New Mexico, Gallup, NM 87301, USA

**Abstract**: In this paper we present theorems and applications of Wallis theorem related to trigonometric integrals.

Let's recall Wallis Theorem:
**Theorem 1.** (Wallis, 1616-1703)
$$\int_0^{\frac{\pi}{2}} \sin^{2n+1} x\,dx = \int_0^{\frac{\pi}{2}} \cos^{2n+1} x\,dx = \frac{2\cdot 4\cdot\ldots\cdot(2n)}{1\cdot 3\cdot\ldots\cdot(2n+1)}.$$

*Proof:* Using the integration by parts, we obtain
$$I_n = \int_0^{\frac{\pi}{2}} \sin^{2n+1} x\,dx = \int_0^{\frac{\pi}{2}} \sin^{2n} x \sin x\,dx = -\cos x \cdot \sin 2nx \Big|_0^{\frac{\pi}{2}} +$$
$$+ 2n\int_0^{\frac{\pi}{2}} \sin^{2n+1} x\left(1 - \sin^2 x\right)dx = 2nI_{n-1} - 2nI_n$$

from where:
$$I_n = \frac{2n}{2n+1} I_{n-1}.$$

By multiplication, we obtain the statement.
We prove in the same manner for $\cos x$.

**Theorem 2.**
$$\int_0^{\frac{\pi}{2}} \sin^{2n} x\,dx = \int_0^{\frac{\pi}{2}} \cos^{2n} x\,dx = \frac{1\cdot 3\cdot\ldots\cdot(2n-1)}{2\cdot 4\cdot\ldots\cdot(2n)} \cdot \frac{\pi}{2}.$$

*Proof:* Same as the first theorem.

**Theorem 3.** If $f(x) = \sum_{k=0}^{\infty} a_{2k} x^{2k}$, then
$$\int_0^{\frac{\pi}{2}} f(\sin x)dx = \int_0^{\frac{\pi}{2}} f(\cos x)dx = \frac{\pi}{2} a_0 + \frac{\pi}{2}\sum_{k=1}^{\infty} a_{2k} \frac{1\cdot 3\cdot\ldots\cdot(2k-1)}{2\cdot 4\cdot\ldots\cdot(2k)}.$$

*Proof:* In the function $f(x) = \sum_{k=0}^{\infty} a_{2k} x^{2k}$ we substitute $x$ by $\sin x$ and then integrate from 0 to $\frac{\pi}{2}$, and we use the second theorem.

**Theorem 4.** If $g(x) = \sum_{k=0}^{\infty} a_{2k+1} x^{2k+1}$, then

$$\int_0^{\frac{\pi}{2}} g(\sin x) dx = \int_0^{\frac{\pi}{2}} g(\cos x) dx = a_1 + \sum_{k=1}^{\infty} a_{2k+1} \frac{2 \cdot 4 \cdot \ldots \cdot (2k)}{1 \cdot 3 \cdot \ldots \cdot (2k+1)}.$$

**Theorem 5.** If $h(x) = \sum_{k=0}^{\infty} a_k x^k$, then

$$\int_0^{\frac{\pi}{2}} h(\sin x) dx = \int_0^{\frac{\pi}{2}} h(\cos x) dx = \frac{\pi}{2} a_0 + a_1 + \sum_{k=1}^{\infty} \left( \frac{\pi}{2} a_{2k} \frac{1 \cdot 3 \cdot \ldots \cdot (2k-1)}{2 \cdot 4 \cdot \ldots \cdot (2k)} + a_{2k+1} \frac{2 \cdot 4 \cdot \ldots \cdot (2k)}{1 \cdot 3 \cdot \ldots \cdot (2k+1)} \right).$$

**Application 1.**

$$\int_0^{\frac{\pi}{2}} \sin(\sin x) dx = \int_0^{\frac{\pi}{2}} \sin(\cos x) dx = \sum_{k=0}^{\infty} (-1)^k \frac{1}{1^2 \cdot 3^2 \cdot \ldots \cdot (2k+1)^2}$$

*Proof:* We use that $\sin x = \sum_{k=0}^{\infty} (-1)^k \frac{x^{2k+1}}{(2k+1)!}$.

**Application 2.**

$$\int_0^{\frac{\pi}{2}} \cos(\sin x) dx = \int_0^{\frac{\pi}{2}} \cos(\cos x) dx = \frac{\pi}{2} \sum_{k=0}^{\infty} \frac{(-1)^k}{4^k (k!)^2}.$$

*Proof:* We use that $\cos x = \sum_{k=0}^{\infty} (-1)^k \frac{x^{2k}}{(2k)!}$.

**Application 3.**

$$\int_0^{\frac{\pi}{2}} sh(\sin x) dx = \int_0^{\frac{\pi}{2}} sh(\cos x) dx = \sum_{k=0}^{\infty} \frac{1}{1^2 \cdot 3^2 \cdot \ldots \cdot (2k+1)^2}.$$

*Proof:* We use that $shx = \sum_{k=0}^{\infty} \frac{x^{2k+1}}{(2k+1)!}$

**Application 4.**

$$\int_0^{\frac{\pi}{2}} ch(\sin x)dx = \int_0^{\frac{\pi}{2}} ch(\cos x)dx = \frac{\pi}{2}\sum_{k=0}^{\infty}\frac{1}{4^k(k!)^2}.$$

*Proof:* We use that $chx = \sum_{k=0}^{\infty}\frac{x^{2k}}{(2k)!}$.

**Application 5.**

$$\sum_{k=1}^{\infty}\frac{1}{k^2} = \frac{\pi^2}{6}$$

*Proof:* In the expression of $\arcsin x = x + \sum_{k=1}^{\infty}\frac{1\cdot 3\cdot ...\cdot(2k-1)x^{2k+1}}{2\cdot 4\cdot ...\cdot(2k)(2k+1)}$ we substitute $x$ by $\sin x$, and use theorem 4. It results that $\frac{\pi^2}{8} = \sum_{k=0}^{\infty}\frac{1}{(2k+1)^2}$.

Because:

$$\sum_{k=1}^{\infty}\frac{1}{k^2} = \sum_{k=0}^{\infty}\frac{1}{(2k+1)^2} + \frac{1}{4}\sum_{k=1}^{\infty}\frac{1}{k^2}$$

we obtain:

$$\sum_{k=1}^{\infty}\frac{1}{k^2} = \frac{\pi}{6}.$$

**Application 6.**

$$\int_0^{\frac{\pi}{2}} \sin x\, ctg(\sin x)dx = \int_0^{\frac{\pi}{2}} \cos x\, ctg(\cos x)dx = \frac{\pi}{2} - \frac{\pi}{2}\sum_{k=1}^{\infty}\frac{B_k}{(k!)^2}$$

where $B_k$ is the k-th Bernoulli type number (see [1]).

*Proof:* We use that $xctgx = 1 - \sum_{k=1}^{\infty}\frac{4^k B_k}{(2k)!}x^{2k}$.

**Application 7.**

$$\int_0^{\frac{\pi}{2}} arctg(\sin x)dx = \int_0^{\frac{\pi}{2}} arctg(\cos x)dx = 1 + \sum_{k=1}^{\infty}(-1)^k \frac{2\cdot 4\cdot ...\cdot(2k)}{1\cdot 3\cdot ...\cdot(2k-1)(2k+1)^2}.$$

*Proof:* We use that $arctgx = \sum_{k=0}^{\infty}(-1)^k\frac{x^{2k+1}}{2k+1}$.

**Application 8.**

$$\int_0^{\frac{\pi}{2}} \arg th(\sin x)dx = \int_0^{\frac{\pi}{2}} \arg th(\cos x)dx = 1 + \sum_{k=1}^{\infty}\frac{2\cdot 4\cdot ...\cdot(2k)}{1\cdot 3\cdot ...\cdot(2k-1)(2k+1)^2}.$$

*Proof:* We use that $\arg th\, x = \sum_{k=0}^{\infty}\frac{x^{2k+1}}{2k+1}$.

**Application 9.**

$$\int_0^{\frac{\pi}{2}} \arg sh(\sin x)dx = \int_0^{\frac{\pi}{2}} \arg sh(\cos x)dx = \sum_{k=1}^{\infty} (-1)^k \frac{1}{(2k+1)^2}.$$

*Proof:* We use that $\arg shx = \sum_{k=0}^{\infty} (-1)^k \frac{1 \cdot 3 \cdot \ldots \cdot (2k-1) x^{2k+1}}{2 \cdot 4 \cdot \ldots \cdot (2k)(2k+1)}$.

**Application 10.**

$$\int_0^{\frac{\pi}{2}} \tg(\sin x)dx = \int_0^{\frac{\pi}{2}} \tg(\cos x)dx = \sum_{k=1}^{\infty} \frac{2^{2k-1}(4^k-1)B_k}{1^2 \cdot 3^2 \cdot \ldots \cdot (2k-1)^2 k}.$$

*Proof:* We use that $\tg x = \sum_{k=1}^{\infty} \frac{2^{2k}(4^k-1)B_k}{(2k)!} x^{2k-1}$.

**Application 11.**

$$\int_0^{\frac{\pi}{2}} \frac{\sin x}{\sin(\sin x)}dx = \int_0^{\frac{\pi}{2}} \frac{\cos x}{\sin(\cos x)}dx = \frac{\pi}{2} + \pi \sum_{k=1}^{\infty} \frac{(2^{2k-1}-1)B_k}{2^{2k}(k!)^2}$$

*Proof:* We use that $\frac{x}{\sin x} = 1 + 2\sum_{k=1}^{\infty} \frac{(2^{2k-1}-1)B_k}{(2k)!} x^{2k}$.

**Application 12.**

$$\int_0^{\frac{\pi}{2}} \frac{\sin x}{sh(\sin x)}dx = \int_0^{\frac{\pi}{2}} \frac{\cos x}{sh(\cos x)}dx = \frac{\pi}{2} + \pi \sum_{k=1}^{\infty} \frac{(2^{2k-1}-1)B_k}{2^{2k}(k!)^2}.$$

*Proof:* We use that $\frac{x}{shx} = 1 + 2\sum_{k=1}^{\infty} (-1)^k \frac{(2^{2k-1}-1)B_k}{(2k)!} x^{2k}$.

**Application 13.**

$$\int_0^{\frac{\pi}{2}} \sec(\sin x)dx = \int_0^{\frac{\pi}{2}} \sec(\cos x)dx = \frac{\pi}{2} + \pi \sum_{k=1}^{\infty} \frac{E_k}{2^{2k+1}(k!)^2},$$

where $E_k$ is the k-th Euler type number (see [1]).

*Proof:* We use that $\sec x = 1 + \sum_{k=1}^{\infty} \frac{E_k}{(2k)!} x^{2k}$

**Application 14.**

$$\int_0^{\frac{\pi}{2}} \sec h(\sin x)dx = \int_0^{\frac{\pi}{2}} \sec h(\cos x)dx = \frac{\pi}{2} + \pi \sum_{k=1}^{\infty} (-1)^k \frac{E_k}{2^{2k+1}(k!)^2}.$$

*Proof:* We use that $\sec h\, x = 1 + \sum_{k=1}^{\infty} (-1)^k \frac{E_k}{(2k)!} x^{2k}$.

**REFERENCES**

[1]   Octav Mayer – Theoria funcțiilor de o variabilă complexă – Ed. Academiei, Bucharest, 1981.
[2]   Mihály Bencze – About Taylor formula – (manuscript).